\documentstyle[amssymb]{amsart}
\title[Macdonald's eigenvalue problem]{Eigenfunctions of  
Macdonald's $q$-difference operator for the root system of type 
$C_n$}
\author{\sc Katsuhisa Mimachi}
\date{}
\numberwithin{equation}{section}
\newtheorem{thm}{Theorem}

\newtheorem{icor}{Corollary}

\newtheorem{lem}{Lemma}
\begin{document}
\maketitle
\begin{abstract}
We construct an integral representation
of eigenfunctions for Macdonald's 
$q$-difference operator associated with 
the root system of type $C_n\,.$   
It is given in terms of a 
restriction of a $q$-Jordan-Pochhammer 
integral.  Choosing a suitable cycle of
the integral, we obtain an integral 
representation of a special case of 
the Macdonald polynomial for 
the root system of type $C_n\,.$ 
\end{abstract}

\section{Introduction}

Macdonald introduced the $q$-difference operators 
\cite{Mac1} to define his orthogonal polynomials
associted with root sytems. 
In the case of a root system of  type $C_n\,,$ his 
$q$-difference operator is given by
$$E=\sum_{a_1,\cdots,a_n=\pm 1}\,
\prod_{1\le i<j\le n}
\frac{1-ty_i^{a_i}y_j^{a_j}}{1-y_i^{a_i}y_j^{a_j}}
\prod_{1\le i\le n}
\frac{1-ty_i^{2a_i}}{1-y_i^{2a_i}}
T_{y_i}^{\frac{a_i}{2}}\,,
$$
where
$$(T_{y_i}f)(y_1,\ldots,y_n)=f(y_1,\ldots,qy_i,\ldots,y_n)\,.$$

The present paper is devoted to study  the eigenvalue 
problem associated with this operator $E\,.$ 
In particular, we construct an integral representation,
which is given by a restriction of a 
$q$-Jordan-Pochhammer integral,  
of eigenfunctions in some special cases. 
It turns out that, taking a 
suitable cycle,  such an integral
expresses the Macdonald polynomial of 
$C_n$ type parametrized by
the partition $(\lambda,0,\ldots,0).$ 
This representaion leads to a more explicit
expression.\smallskip 

Here we recall the definition of the Macdonald polynomial
$P_\mu(y|q,t)$ associated with
the root system of type $C_n\,.$ It is the eigenfunction of $E$ with
respect to the eigenvalue
$$c_{\mu}=
q^{-\frac{1}{2}(\mu_1+\cdots+\mu_n)}
\prod_{i=1}^n (1+q^{\mu_i}\,t^{n-i+1})$$
of the form
$$P_\mu(y|q,t)=m_\mu+\sum_{\nu<\mu}a_{\mu\;\nu }m_\nu\,,$$
where $\mu=(\mu_1,\ldots,\,\mu_n)$ is a partition, 
a sequence of non-negative integers in decreasing order, 
$m_\mu=\sum_{\nu\in W(C_n)\mu}e^\nu$ with $W(C_n)$ the Weyl
group of type $C_n$  and $\nu<\mu$ is defined to be 
$\mu-\nu\in Q^+$ with $Q^+$ the positive cone of the root 
lattice.\medskip 


Besides the $A_{n-1}$ case, the solution of the eigenvalue problem for
the Macdonald operator is not well studied (See \cite{Mi1,Mi2,MN}
and \cite{Mac3}).
We expect that this paper represents
a step toward understanding the $BC_n$ type 
Macdonald polynomials\cite{Koorn}. 
It is noteworthy that even in 
the classical ($q$=1) case was not previously 
known that such an integral gives spherical 
functions associated with the root system $C_n.$
For related works on $BC_n$ type spherical 
functions, we refer the reader to \cite{DG} 
and references therein.\medskip

Throughout this paper, $q$ is regarded as a real 
number satisfying $0< q <1\,,$ and $t=q^k$ where
$k\in {\Bbb Z}_{\ge 1}\,.$ 

\section{A restriction of 
a $q$-Jordan-Pochhammer integral}

Let us introduce a 1-form 
\begin{equation}
\begin{split}
\Phi&=x^{\lambda}\prod_{1\le j\le n}
\frac{ (ty_j/x;q)_\infty\; (ty_j^{-1}/x;q)_\infty}
{(\;y_j/x;q)_\infty\; (\;y_j^{-1}/x;q)_\infty}\frac{dx}{x}\\
&=x^{\lambda}\prod_{1\le j\le n}
\frac{1}
{(y_j/x;q)_k\;(y_j^{-1}/x;q)_k}\frac{dx}{x}\,,
\end{split}
\end{equation}
where $\lambda\in {\Bbb Z}_{\ge 0}, 
(a;q)_\infty=\prod_{i\ge 0}(1-aq^i)\,$ and 
$(a;q)_m=(a;q)_\infty/(q^ma;q)_\infty\,.$
This can be regarded as a 1-form corresponding to a restriction of 
a $q$-Jordan-Pochhammer integral 
$$x^{\lambda}\prod_{1\le j\le 2n}
\frac{ (ty_j/x;q)_\infty }
{(y_j/x;q)_\infty }\frac{dx}{x}\,,$$
which is studied in \cite{Mi0} and \cite{AKM}.\medskip

Our first result is the following:

\begin{thm} For any cycle ${\cal C}\,,$  
the function $\int_{\cal C}\Phi$ satisfies the equation
$$E\int_{\cal C}\Phi=c_{(\lambda,0,\ldots,0)}\int_{\cal C}\Phi\,.$$
\end{thm}
This implies that linearly independent solutions are obtained 
by choosing several cycles. Indeed, if we put 
$C_i^{(+)}$ (or $C_i^{(-)}$) for each $i=1,\ldots,n$
to be a path with the counterclockwise direction so that
the poles at $w=y_i, y_iq,\ldots, y_iq^{k-1}$ 
(or $w=y_i^{-1},y_i^{- 1}q,\ldots, y_i^{- 1}q^{k-1}$, respectively) 
are inside the path and other poles from $\Phi$ are outside, 
we have the following rational solutions: 
{\allowdisplaybreaks
\begin{equation}
\begin{split}
\frac{1}{2\pi\sqrt{-1}}\int_{C_i^{(+)}}\Phi&=y_i^{\lambda}
\frac{1}{\displaystyle (q;q)_{k-1}
\prod\begin{Sb}
1\le j\le n\\
j\ne i
\end{Sb}
(y_j/y_i;q)_k
\prod\begin{Sb}
1\le j\le n
\end{Sb}
(y_j^{-1}y_i^{-1};q)_k}\\
&\quad\times
\sum_{l=0}^{k-1}\prod_{j=1}^n
\frac{\displaystyle 
(t^{-1}qy_i/y_j;q)_l\;(t^{-1}qy_iy_j;q)_l}
{\displaystyle  
(qy_i/y_j;q)_l\;(qy_iy_j;q)_l}
(t^{2n}q^\lambda)^l\\
\end{split}
\end{equation}}
and
{\allowdisplaybreaks
\begin{equation}
\begin{split}
\frac{1}{2\pi\sqrt{-1}}\int_{C_i^{(-)}}\Phi&=y_i^{-\lambda}
\frac{1}{\displaystyle (q;q)_{k-1}
\prod\begin{Sb}
1\le j\le n\\
j\ne i
\end{Sb}
(y_i/y_j;q)_k
\prod\begin{Sb}
1\le j\le n
\end{Sb}
(y_jy_i;q)_k}\\
&\quad\times
\sum_{l=0}^{k-1}\prod_{j=1}^n
\frac{\displaystyle 
(t^{-1}qy_j/y_i;q)_l\;(t^{-1}qy_i^{-1}y_j^{-1};q)_l}
{\displaystyle  
(qy_j/y_i;q)_l\;(qy_i^{-1}y_j^{-1};q)_l}
(t^{2n}q^\lambda)^l\,.\\
\end{split}
\end{equation}}\noindent
The calculation is carried out by means of the residue calculus.\medskip

Since $\lambda$ is a non-negative integer, the sum of the pathes
$\sum_{i=1}^nC_i^{(+)}+\sum_{i=1}^nC_i^{(-)}$ is homologous to a path
$C$ which circles the origin in the positive sense so that all poles from 
$\Phi$ are inside the path. The integral on this cycle $C$ gives the
Macdonald polynomial $P_{(\lambda,0,\ldots,0)}(y|q,t)\,.$
\begin{thm} If the cycle $C$ is that above, we have
\begin{equation}
\frac{1}{2\pi\sqrt{-1}}\int_C\Phi=
\frac{(t;q)_\lambda}{(q;q)_\lambda}
P_{(\lambda,0,\ldots,0)}(y|q,t).
\end{equation}
\end{thm}\medskip
Moreover, applying the $q$-binomial theorem 
$$\sum_{i\ge 0}\frac{(a;q)_i}{(q;q)_i}z^i
=\frac{(az;q)_\infty}{(z;q)_\infty}\quad(|z|<1)
$$
with the residue calculus to our integral, we obtain 
an exact expression of  
$P_{(\lambda,0,\ldots,0)}(y|q,t)\,.$\par\smallskip

\begin{icor}
$$P_{(\lambda,0,\ldots,0)}(y|q,t)
=\frac{(q;q)_\lambda}{(t;q)_\lambda}\sum
\begin{Sb}
i_1+\cdots+i_{2n}=\lambda\\[3pt]
i_1,\,\ldots,\,i_{2n}\ge 0
\end{Sb}
\frac{(t;q)_{i_1}\cdots\,(t;q)_{i_{2n}}}
{(q;q)_{i_1}\cdots\,(q;q)_{i_{2n}}}
y_1^{i_1-i_{2n}}y_2^{i_2-i_{2n-1}}\cdots\,y_n^{i_n-i_{n+1}}\,.
$$ 
\end{icor}

\section{Proof of Theorem 1}
\begin{lem} We have
$$\sum_{a_1,\ldots,a_n=\pm 1}\;\prod_{1\le i<j\le n}
\frac{1-ty_i^{a_i}y_j^{a_j}}{1-y_i^{a_i}y_j^{a_j}}
\prod_{1\le i\le n}
\frac{1-ty_i^{2a_i}}{1-y_i^{2a_i}}=
\prod_{i=1}^n(1+t^i)\,.$$
\end{lem}
From the formula by Macdonald\cite{Mac0} about the 
Poincar\'{e} series of Coxter systems, we have
\begin{equation}
\sum_{w\in W(C_n)}\;w\left\{\prod_{1\le i<j\le n}
\frac{1-ty_iy_j}{1-y_iy_j}
\prod_{1\le i\le n}
\frac{1-ty_i}{1-y_i}\right\}=
\frac{\prod_{i=1}^n(1-t^{2i})}{(1-t)^n}
\end{equation}
and
\begin{equation}
\sum_{w\in W(A_{n-1})}\;w\left\{\prod_{1\le i<j\le n}
\frac{1-ty_i/y_j}{1-y_i/y_j}\right\}=
\frac{\prod_{i=1}^n (1-t^{i})}{(1-t)^n}\,.
\end{equation}
Here $W(C_n)$ or $W(A_{n-1})$ denotes the Weyl group of
the root system of type $C_n$ or $A_{n-1}\,,$ respectively.
By applying the formula (3.2) to (3.1), we obtain
{\allowdisplaybreaks
\begin{align*}
&\sum_{w\in W(C_n)}\;w\left\{\prod_{1\le i<j\le n}
\frac{1-ty_iy_j}{1-y_iy_j}
\prod_{1\le i\le n}
\frac{1-ty_i}{1-y_i}\right\} \\
&=\frac{\prod_{i=1}^n (1-t^{i})}{(1-t)^n}
\sum_{a_1,\ldots,a_n=\pm 1}\;\prod_{1\le i<j\le n}
\frac{1-ty_i^{a_i}y_j^{a_j}}{1-y_i^{a_i}y_j^{a_j}}
\prod_{1\le i\le n}
\frac{1-ty_i^{2a_i}}{1-y_i^{2a_i}}\,.
\end{align*}
}
Hence we derive the desired relation.\quad\qed\medskip
\begin{lem}
\begin{align}
&\sum_{a_1,\ldots,a_n=\pm 1}\;\prod_{1\le i<j\le n}
\frac{1-ty_i^{a_i}y_j^{a_j}}{1-y_i^{a_i}y_j^{a_j}}
\prod_{1\le i\le n}
\frac{(1-ty_i^{2a_i})(1-y_i^{a_i}/x)}
{(1-y_i^{2a_i})(1-ty_i^{a_i}/x)}\\\notag
&\quad=\prod_{i=1}^{n-1}(1+t^i)\;
\left\{1+t^n\prod_{i=1}^n
\frac{(1-y_i/x)(1-y_i^{-1}/x)}
{(1-ty_i/x)(1-ty_i^{-1}/x)}\right\}\,.
\end{align}
\end{lem}
{\it Proof.}
We prove the desired equality by means of partial 
fraction decompositions.
Firstly, let us take the residue of the left-hand side of (3.3) 
at $x=ty_1\,:$
\begin{align*}
&\text{Res}_
{\,x=ty_1}
\Biggl\{
\sum_{a_1,\ldots,a_n=\pm 1}\;\prod_{1\le i<j\le n}
\frac{1-ty_i^{a_i}y_j^{a_j}}{1-y_i^{a_i}y_j^{a_j}}
\prod_{1\le i\le n}
\frac{(1-ty_i^{2a_i})(1-y_i^{a_i}/x)}
{(1-y_i^{2a_i})(1-ty_i^{a_i}/x)}\Biggr\}\frac{dx}{x}\\
&=\frac{(1-t^{-1})(1-ty_1^2)}{1-y_1^2}
\sum_{a_2,\ldots,a_n=\pm 1}\;
\prod_{2\le j\le n}
\frac{1-ty_1y_j^{a_j}}{1-y_1y_j^{a_j}}\;\\
&\qquad\times\prod_{2\le i<j\le n}
\frac{1-ty_i^{a_i}y_j^{a_j}}{1-y_i^{a_i}y_j^{a_j}}
\prod_{2\le i\le n}
\frac{(1-ty_i^{2a_i})(1-t^{-1}y_i^{a_i}/y_1)}
{(1-y_i^{2a_i})(1-y_i^{a_i}/y_1)}\\
&=t^{1-n}\frac{(1-t^{-1})(1-ty_1^2)}{1-y_1^2}
\prod_{2\le j\le n}
\frac{(1-ty_1y_j)(1-ty_1/y_j)}
{(1-\;y_1y_j)(1-\;y_1/y_j)}\\
&\qquad\times\sum_{a_2,\ldots,a_n=\pm 1}\;\prod_{2\le i<j\le n}
\frac{1-ty_i^{a_i}y_j^{a_j}}{1-\;y_i^{a_i}y_j^{a_j}}
\prod_{2\le i\le n}
\frac{1-ty_i^{2a_i}}
{1-\;y_i^{2a_i}}\,.
\end{align*}
This is equal to 
\begin{equation}
t^{1-n}\frac{(1-t^{-1})(1-ty_1^2)}{1-y_1^2}
\prod_{i=1}^{n-1}(1+t^i)
\prod_{2\le j\le n}
\frac{(1-ty_1y_j)(1-ty_1/y_j)}
{(1-\;y_1y_j)(1-\;y_1/y_j)}\,,
\end{equation}
from Lemma 1.\medskip

Secondly, by noticing that
\begin{align*}
&\text{Res}_
{\,x=ty_1}
\prod_{1\le i\le n}
\frac{(1-\;y_i/x)(1-\;y_i^{-1}/x)}{(1-ty_i/x)(1-ty_i^{-1}/x)}
\frac{dx}{x}\\
&=t^{1-2n}\frac{(1-t^{-1})(1-ty_1^2)}{1-y_1^2}
\prod_{2\le j\le n}
\frac{(1-ty_1y_j)(1-ty_1/y_j)}
{(1-\;y_1y_j)(1-\;y_1/y_j)}\,,
\end{align*}
we know that the residue of the right-hand side of (3.3) at $x=ty_1$ 
is equal to (3.4).
Hence, the symmetry of (3.3) with respect to
the variables $y_1^{\pm 1},\ldots,y_n^{\pm 1}\,$ 
leads to the fact that
the residues of both sides of (3.3) at each
$x=ty_i\quad (i=1,\ldots,n)$ or 
$x=ty_i^{-1}\quad (i=1,\ldots,n)\,$ are equal.\medskip

On the other hand, if $x$ goes to $\infty$,  
the left-hand side of (3.3) tends to
$$
\sum_{a_1,\ldots,a_n=\pm 1}\;\prod_{1\le i<j\le n}
\frac{1-ty_i^{a_i}y_j^{a_j}}{1-y_i^{a_i}y_j^{a_j}}
\prod_{1\le i\le n}
\frac{1-ty_i^{2a_i}}
{1-\;y_i^{2a_i}}\,,
$$
which is equal to $\prod_{i=1}^n(1+t^i)$ by Lemma 1, 
and the right-hand side of (3.3) tends also to
$\prod_{i=1}^n(1+t^i)\,.$
This completes the proof of Lemma 2.\quad\qed\medskip

Let us proceed to  prove our Theorem 1.

Note that
{\allowdisplaybreaks
\begin{equation*}
\begin{split}
&\prod_{1\le i\le n}
T^{\frac{a_i}{2}}_{y_i}
\int\Phi=\int_C\,x^\lambda\prod_{i=1}^n
\frac{\displaystyle
\biggl(q^{\frac{a_i}{2}}\frac{ty_i}{x}\biggr)_\infty
\biggl(q^{-\frac{a_i}{2}}\frac{ty_i^{-1}}{x}\biggr)_\infty}
{\displaystyle
\biggl(q^{\frac{a_i}{2}}\frac{\,y_i}{x}\biggr)_\infty
\biggl(q^{-\frac{a_i}{2}}\frac{\,y_i^{-1}}{x}\biggr)_\infty}
\frac{dx}{x}\\
&=q^{-\frac{\lambda}{2}}
\int_C\,x^\lambda\prod_{i=1}^n
\frac{\displaystyle
\biggl(q^{\frac{1}{2}(1+a_i)}\frac{ty_i}{x}\biggr)_\infty
\biggl(q^{\frac{1}{2}(1-a_i)}\frac{ty_i^{-1}}{x}\biggr)_\infty}
{\displaystyle
\biggl(q^{\frac{1}{2}(1+a_i)}\frac{\,y_i}{x}\biggr)_\infty
\biggl(q^{\frac{1}{2}(1-a_i)}\frac{\,y_i^{-1}}{x}\biggr)_\infty}
\frac{dx}{x}\\
&=q^{-\frac{\lambda}{2}}
\int_C\,\prod_{i=1}^n
\frac
{\displaystyle 1-y_i^{a_i}/x}
{\displaystyle 1-ty_i^{a_i}/x}
\Phi\,,
\end{split}
\end{equation*}}
where the second equality is given by the change of 
integration variable such that 
$x\mapsto q^{-\frac{1}{2}}x\,.$\medskip

Therefore, by using Lemma 2, we obtain
{\allowdisplaybreaks
\begin{equation*}
\begin{split}
&E\int_C\Phi=q^{-\frac{\lambda}{2}}\int_C\Biggl\{
\sum_{a_1,\ldots,a_n=\pm 1}\;\prod_{1\le i<j\le n}
\frac{1-ty_i^{a_i}y_j^{a_j}}{1-y_i^{a_i}y_j^{a_j}}
\prod_{1\le i\le n}
\frac{(1-ty_i^{2a_i})(1-y_i^{a_i}/x)}
{(1-y_i^{2a_i})(1-ty_i^{a_i}/x)}\Biggr\}\Phi\\
&=q^{-\frac{\lambda}{2}}
\prod_{i=1}^{n-1}(1+t^i)\int_C
\left\{1+t^n\prod_{i=1}^n
\frac{(1-y_i/x)(1-y_i^{-1}/x)}
{(1-ty_i/x)(1-ty_i^{-1}/x)}\right\}\Phi\\
&= 
(q^{-\frac{\lambda}{2}}+q^{\frac{\lambda}{2}}t^n)\,
\prod_{i=1}^{n-1}(1+t^i)\,
\int_C\Phi\,.
\end{split}
\end{equation*}}\noindent
Here, to derive the third equality, we have used 
the relation
$$\int_C\Phi=q^{-\lambda}
\int_C
\prod_{i=1}^n
\frac{(1-y_i/x)(1-y_i^{-1}/x)}
{(1-ty_i/x)(1-ty_i^{-1}/x)}\Phi\,,$$
which is given by the change of integration variable
such that $x\mapsto q^{-1}x\,.$\par
This completes the proof of Theorem 1.

\bigskip
\begin{flushleft}
\begin{sc}
Katsuhisa Mimachi\\
Department of Mathematics\\
Kyushu University 33\\
Hakozaki, Fukuoka 812-81\\
Japan\\
\end{sc}
\smallskip 
{\it E-mail address\/}: mimachi{\char'100}
math.kyushu-u.ac.jp\\
\end{flushleft}

\begin{thebibliography}{99}
\bibitem{AKM} Aomoto, K., Kato, Y., Mimachi, K.: 
A solution of the Yang-Baxter equation as connection
coefficients of a holonomic $q$-difference system. 
Internat. Math. Res. Notices 
{\bf 1992}, No.1 , 7-15 
\bibitem{DG} Debiard, A. and Gaveau, B.: Integral formulas for the spherical 
polynomials of a root system of type $BC_2$. Jour. of Funct. Anal. {\bf 119}, 
401-454 (1994)
\bibitem{Koorn} Koornwinder, T. H.: {\it Askey-Wilson polynomials 
for root systems of type BC,} in
Hypergeometric functions on domains of positivity, Jack polynomials 
and applications,
D.St.P.Richards(ed.), Contemp. Math. {\bf 138}, Amer. Math. Soc. 
(1992), pp.189-204 
\bibitem{Mac0} Macdonald, I.G.: The Poincar\'e series of a Coxeter group,
Math. Ann. {\bf 199}, 161-174 (1972)
\bibitem{Mac1} Macdonald, I.G.: {\em A new class of symmetric 
functions,} in Actes S\'eminaire Lotharingen,
Publ. Inst. Rech. Math. Adv., Strasbourg, 1988, 131-171
\bibitem{Mac2} Macdonald, I.G.: Affine Hecke algebras and 
orthogonal polynomials. S\'eminaire
BOURBAKI, 47\`eme ann\'ee, 1994-95,  $n^{\circ} 797$
\bibitem{Mac3} Macdonald, I.G.: 
{\em Symmetric Functions and Hall Polynomials} (Second Edition),
Oxford Mathematical Monographs, Clarendon Press, Oxford, 1995.
\bibitem{Mi0} Mimachi, K.: Connection problem in holonomic 
$q$-difference system associated with
a Jackson integral of Jordan-Pochhammer type. Nagoya Math. J. 
{\bf 116}, 149-161 (1989)
\bibitem{Mi1} Mimachi, K.:
A solution to quantum Knizhnik-Zamolodchikov
equations and its application to eigenvalue problems of 
the Macdonald type,
Duke Math.\,J.\, {\bf 85}, 635-658 (1996).
\bibitem{Mi2} Mimachi, K.:
Rational solutions to eigenvalue problems 
of the Macdonald type, in preparation.
\bibitem{MN} Mimachi, K. and Noumi, M.: An integral representation
of eigenfunctions for Macdonald's $q$-difference operators, 
T\^ohoku Math.\,J.\, {\bf 49}, 517-525 (1997).
\end{thebibliography}
\end{document}